\documentclass{amsart}

\usepackage{amssymb}

\theoremstyle{plain}

\numberwithin{equation}{section}

\newtheorem*{teo}{Theorem A}
\newtheorem*{teob}{Theorem B}

\newtheorem{hypo}[equation]{Hypotheses}

\newtheorem{lem}[equation]{Lemma}

\newtheorem{prop}[equation]{Proposition}

\newtheorem{claim}[equation]{Claim}

\newcommand{\modu}{\operatorname{mod}}

\newcommand{\Aff}{\operatorname{Aff}}

\theoremstyle{definition}

\begin{document}	
\title{On nilpotent groups and conjugacy classes}

\author{Edith Adan-Bante}

\address{University of Southern Mississippi Gulf Coast, 730 East Beach Boulevard,
 Long Beach MS 39560}

\email{Edith.Bante@usm.edu}


\subjclass{20d15}

\date{2005}
\begin{abstract} Let $G$ be a nilpotent group 
and $a\in G$. Let 
$a^G=\{g^{-1}ag\mid g\in G\}$ be the conjugacy class of $a$ in $G$.
Assume that $a^G$ and $b^G$ are conjugacy classes of $G$ 
with the property that $|a^G|=|b^G|=p$, where $p$ is 
an odd prime number.  
 Set $a^G b^G=\{xy\mid x\in a^G, y\in b^G\}$.  Then 
either $a^G b^G=(ab)^G$ or $a^G b^G$ is the union of at least $\frac{p+1}{2}$ distinct
conjugacy classes. As an application of the previous 
result, given any nilpotent group $G$
and any conjugacy class $a^G$ of size $p$,  we describe  the square
$a^G a^G$ of $a^G$ in terms of conjugacy classes of $G$.
\end{abstract}
\maketitle

\begin{section}{Introduction}

Let $G$ be a finite group, $a\in G$  and $a^G=\{a^g\mid g\in G\}$ 
be the conjugacy class of $a$ in $G$.
Let $X$ be a $G$-invariant subset of $G$, i.e. 
$X^g=\{x^g\mid x\in X\}=X$ for all $g\in G$.
  Then $X$ can be expressed as a union of 
  $n$ distinct conjugacy classes of $G$, for some integer $n>0$. Set
 $\eta(X)=n$.
 
  We can check that given any two conjugacy classes
 $a^G$ and $b^G$ of $G$, the product $a^G b^G=\{a^g b^h\mid g, h\in G\}$ is a 
 $G$-invariant set. Thus $\eta(a^G b^G)$ is the number of distinct conjugacy 
 classes of $G$ such that $a^G b^G$ is the union of those classes.
 In this 
 note we  study the product $a^G b^G$ when $G$ is a nilpotent group,
  and $|a^G|=|b^G|=p$ for some prime $p$.
  The main result is the following
  
\begin{teo} Let $p$ be an odd prime number.
Let $G$ be a finite nilpotent group, $a^G$ and
 $b^G$ be conjugacy classes with the property 
that $|a^G|=|b^G|=p$.
 Then either $a^G b^G =(ab)^G$
 or $a^G b^G$ is the union of 
 at least $\frac{p+1}{2}$ distinct conjugacy classes.
\end{teo}

For a fix prime $p>3$, we want to emphasize that Theorem A is a statement about 
 a ``gap" among the possible values that $\eta(a^G b^G)$ can take for any
finite nilpotent group $G$ and any conjugacy classes $a^G$ and $b^G$ of size $p$. 
We wonder if there exist any other  ``gaps".
 
  Given any prime $p>2$ and any integer $n>0$,
we provide an example in Proposition \ref{basico}
of a $p$-group $G$ and a conjugacy class $a^G$ of $G$ 
such that $|a^G|=p^n$ and 
$\eta(a^Ga^G)=\frac{p+1}{2}$. Thus the bound of Theorem A is optimal.

  In Proposition \ref{sizetwo} it is proved that if 
  $G$ is a finite group, $a^G$ and $b^G$ are conjugacy classes of size two, then
  $\eta(a^G b^G)=1$ or $\eta(a^G b^G)=2$. 
  
An application of Theorem A is the following

\begin{teob} 
Let $p$ be an odd prime number.
Let $G$  be a nilpotent group and $a^G$ be a conjugacy class of size $p$. 
Then one of the following holds

{\bf i)}  $a^G a^G=(a^2)^G$. Therefore $[a,G]=[a^2, G]$ and $[a, G]$ 
 is a normal subgroup of $G$.

{\bf ii)} $a^G a^G$ is the union of exactly $\frac{p+1}{2}$ distinct conjugacy classes
of $G$ of size $p$. 
 \end{teob}

 Given any prime $p>2$, we provide
  an example in Proposition \ref{super}
  of a finite
 supersolvable group $G$ and a conjugacy class $a^G$ of $G$ of size
 $p$ with $\eta(a^Ga^G)=2$. Thus  Theorem B
may not remain true if $G$ is only a  
supersolvable group.  

For similar results in character theory, see \cite{edithsquare}.

 {\bf Acknowledgment.} I would like to thank Professor Everett C. Dade  
 for very useful emails, and for helping me to improve the presentation of 
 this note.  
 
 \end{section}

\begin{section}{Proof of Theorem A}   
 Denote by 
${\bf C}_G(a)=\{g\in G\mid a^g=a\}$ the centralizer of $a$ in $G$.
Also let $[a,g]=a^{-1} a^g$ and $[a,G]=\{[a,g]\mid g\in G\}$. Since $a^g=a[a,g]$ for
any $g\in G$, we have that 
$a^G=\{ax\mid x\in [a,G]\}$. In particular, we get that $|a^G|=|[a,G]|$.

Through this note, we will use the well know fact that $|a^G|=|G:{\bf C}_G(a)|$. 

\begin{lem}\label{counting} Let $p$ be an odd prime number, and $r,s$ and 
$t$ be integer numbers.
The set $\{ri^2+si+t\modu p\mid i=0, \ldots ,p-1\}$ has either one element
or at least $\frac{p+1}{2}$ elements.
\end{lem}
\begin{proof}
If $r\equiv 0\modu p$, then $|\{ri^2+si+t\modu p\mid i=0, \ldots ,p-1\}|=|\{si+t\modu p\mid i=0,\ldots,p-1\}|$. If $s\equiv 0\modu p$, in addition, then  $|\{si+t \modu p\mid i=0,\ldots,p-1\}|=|\{t\}|=1$.
Otherwise $|\{si+t\modu p\mid i=0,\ldots,p-1\}|=p$. Thus we may assume that 
$r\not\equiv 0\modu p$. 

If $ri^2+si+t\equiv rj^2+sj+t\modu p$, then $r(i^2-j^2)\equiv r(i-j)(i+j)\equiv -s(i-j)\modu p$. Therefore either
$i-j\equiv 0\modu p$ or $r(i+j)\equiv -s\modu p$. Thus either 
$i\equiv j \modu p$ or $i\equiv -j - \frac{s}{r} \modu p$.
 We conclude that the set $\{ri^2+si+t\modu p\mid i=0, \ldots ,p-1\}$ either has only one
 element or it
has at least $\frac{p}{2}$ elements, 
and therefore at least  $\frac{p+1}{2}$ elements.
\end{proof}

\begin{lem}\label{lemma1}
Let $G$ be a finite $p$-group
 and $N$ be a normal subgroup of $G$. 
Let $a$ and $b$ be elements
of $G$. 
If $(aN)^{G/N}\cap (bN)^{G/N}= \emptyset$
 then $a^G\cap b^G= \emptyset$.
 Thus $\eta((aN)^{G/N} (bN)^{G/N})\leq 
  \eta(a^G b^G )$.
\end{lem}
\begin{proof}
See Lemma 2.1 of  \cite{edith2}.
\end{proof}

\begin{lem}\label{productsn} Let $G$ be a finite group and  $a,b\in G$. Then
\begin{equation*}\label{asquare}
a^G b^G =ab[a^b,G][b,G].
\end{equation*}
\end{lem}
\begin{proof}
See Lemma 2.1 of \cite{edithhomogen}.
\end{proof}

\begin{lem}\label{samesize}
Let $G$ be a finite group, $a^G$ and $b^G$ be  conjugacy classes of $\,G$. Assume that 
$|a^G|=|b^G|=|(ab)^G|$.
 Then  $a^G b^G=(ab)^G$ if and only if  $[ab,G]=[a,G]=[b,G]$ and $[ab,G]$ is a normal 
subgroup of $\,G$. 
   \end{lem}
\begin{proof}
Assume that $a^G b^G=(ab)^G$.
By Lemma \ref{productsn} we have that $[a^b, G][b,G]=[ab,G]$. Observe that
$1_G\in [a^b, G]$  and $1_G\in[b,G]$. Observe also that
$|[a^b, G]|=|[a,G]^{b}|=|[a,G]|$.  Since $|a^G|=|b^G|=|(ab)^G|$,
and $[a^b, G][b,G]=[ab,G]$, it follows that $[ab,G]=[a^b,G]=[b,G]$.

Let $u,v\in [ab,G]$. 
Since $[ab,G]=[a^b,G]=[b,G]$, we have that
$u\in [a^b,G]$ and  $v\in [b,G]$, and therefore
 $uv\in [a^b, G][b,G]=[ab,G]$.
Since $1_G\in [ab,G]$ and $G$ is a finite group, it follows that
$[ab,G]$ is a subgroup of $G$. By Lemma 2.2 of \cite{edithhomogen} we have that
 $[ab,G]$ is a normal subgroup of $G$ since it is a subgroup of $G$. Thus $[a,G]= ([a,G]^{b})^{b^{-1}}=[a^b, G]^{b^{-1}}=[ab,G]^{b^{-1}}=[ab,G]$ and therefore
 $[ab,G]=[a,G]=[b,G]$ is a normal subgroup of $G$.
 
 If $[ab,G]=[a,G]=[b,G]$ and $[ab,G]$ is a normal 
subgroup of $\,G$, then 
$[a^b,G]=[a,G]^{b}=[a,G]$ and so 
$ [a^b,G][b,G]=[ab,G]$. By Lemma \ref{productsn} it follows that
$a^G b^G=(ab)^G$ and the proof is complete.
\end{proof}

\begin{lem}\label{centralequal}
   Let $G$ be a finite group, $a^G$ and $b^G$ be conjugacy classes of $\,G$.
 If $ a^G b^G\cap{\bf Z}(G)\neq \emptyset $ then 
  ${\bf C}_G(b)={\bf C}_G(a^g)$ for some $g\in G$.
  If, in addition,  ${\bf C}_G(a)$ is a normal subgroup of $G$ then 
  ${\bf C}_G(b)={\bf C}_G(a)$.
 \end{lem}
 \begin{proof}
Let $z\in a^G b^G\cap{\bf Z}(G)$. Then $z=a^g b$ for some $g\in G$ and therefore
$a^g=b^{-1}z$. Observe that 
 \begin{eqnarray*}
   {\bf C}_G(a^g)&=&\{h\in G\mid (a^g)^h=a^g\}=\{h\in G\mid (b^{-1}z)^h=b^{-1}z\}\\
    &=&\{h\in G\mid (b^{-1})^h z=b^{-1}z\}=\{h\in G\mid (b^{-1})^h=b^{-1}\}\\
    &=&\{h\in G\mid b^h=b\}={\bf C}_G(b).
    \end{eqnarray*}
    
   We can check that ${\bf C}_G(a^g)=({\bf C}_G(a))^g$ and the second statement 
   follows.
   \end{proof}
   \begin{lem}\label{squaresofsize2}
 Let $G$ be a finite group, $a^G$ and $b^G$ be conjugacy classes of $\,G$ with 
 ${\bf C}_G(a)={\bf C}_G(b)$
 and
 $|a^G|=2$.
Then $\eta(a^G b^G)= 2$.
\end{lem}
\begin{proof}
See Proposition 3.4 of \cite{edithhomogen}.
\end{proof}
\begin{prop}\label{sizetwo}
Let $G$ be a finite group, $a^G$ and $b^G$ be conjugacy classes of $\,G$ with 
$|a^G|=|b^G|=2$. Then either $a^G b^G=(ab)^G$ or $a^G b^G$ is the union
of two distinct conjugacy classes of $G$. 
\end{prop}
\begin{proof}
Since $|G:{\bf C}_G(a)|=|a^G|=2$, we have that ${\bf C}_G(a)\triangleleft G$. 
  If   $a^G b^G\cap {\bf Z}(G)\neq \emptyset$ then by Lemma \ref{centralequal}
  we have that ${\bf C}_G(a)={\bf C}_G(b)$. Thus by Lemma \ref{squaresofsize2} we have that
  $\eta(a^G b^G)=2$. We may assume then 
  that $a^G b^G\cap {\bf Z}(G)=\emptyset$ and so if $c^G\subseteq a^G b^G$, then 
  $|c^G|\geq 2$. Since $|a^G b^G|\leq 4$, it follows then that 
  $\eta(a^G b^G)\leq 2$. 
 \end{proof} 

\begin{lem}\label{intersection}
Let $G$ be a finite group, $b^G$ be a conjugacy class of $\,G$ 
and $N\subseteq {\bf Z}(G)$. Then $b^G n$ is a conjugacy class of $G$ for any 
$n\in N$. Also given any $n_1$ and $n_2$ in $N$, $b^G n_1=b^Gn_2$ if and only if 
$n_1 n_2^{-1}\in [b,G]$. Therefore, if $N\cap [b,G]=\{1_G\}$, then 
 $b^G n_1=b^Gn_2$ if and only if $n_1=n_2$.
\end{lem}
\begin{proof}
Since $N\subseteq {\bf Z}(G)$, $b^G n=(bn)^G$. Given any $n_1, n_2$ in $N$, observe 
 that  $b^G n_1=b^Gn_2$ if and only if
$b^G= b^G n_1 n_2^{-1}$. Also $b^G= b^G n_1 n_2^{-1}$ if and only if there exist 
some $g\in G$ such that 
 $b^g=bn_1 n_2^{-1}$, i.e. if and only if $n_1n_2^{-1}\in [b,G]$.
\end{proof}
\begin{lem}\label{observation1}
Let $K$ and $L$ be finite groups, $X$ be a $K$-invariant subset of $K$ and
$Y$ be a $L$-invariant subset of $\,L$. Let $G=K\times L$ be the direct product of 
$K$ and $L$. Set
$X\times Y=\{(x,y)\in G\mid x\in X, y\in Y\}$.  Then $X\times Y$ is a $G$-invariant
subset of $G$ and 
$\eta(X\times Y)=\eta(X) \eta(Y)$. 
In particular, given any conjugacy classes $a^K$ and $b^K$ of $\,K$, any 
conjugacy classes $c^L$ and $d^L$ of $\,L$, we have that 
$\eta((a,c)^G(b,d)^G)=\eta(a^K b^K)\eta(c^L d^L)$.
\end{lem}
\begin{proof}
Let $X=\cup_{i=1}^{\eta(X)} x_i^K$ and $Y=\cup_{j=1}^{\eta(Y)} y_j^L$.
Since $G=K\times L$, given any $i\in \{1, \ldots, \eta(X)\}$
and any $j\in\{1, \ldots, \eta(Y)\}$,
we can check that 
 $$(x_i, y_j)^G =\{(x,y)\mid x\in x_i^K, y\in y_j^L\}=(x_i^K, y_j^L).$$
Observe that 
  $(x_i^K, y_j^L)=(x_{i_o}^K, y_{j_0}^L)$ if and only if
  $x_i^K=x_{i_0}^K$ and $y_j^L = y_{j_0}^L$. Since
  $X\times Y$ is the 
  union of $(x_i^K, y_j^L)$ for $i=1, \ldots, \eta(X)$, and  $j=1, \ldots, \eta(Y)$, we 
  conclude that
  $\eta(X\times Y)=\eta(X) \eta(Y)$. 
\end{proof}

\begin{proof}[Proof of Theorem A]
Since $p$ is a prime number and $G$ is a nilpotent group, by Lemma \ref{observation1}
 we may assume that $G$ is a $p$-group.

Since $|a^G|=p\,$ and $|G:{\bf C}(a)|=|a^G|$, 
we have  that ${\bf C}_G(a)$ is a normal subgroup of $G$ of index $p$. Similarly,
since $|b^G|=p$, we have that ${\bf C}_G(b)$ is a normal subgroup of $G$ of index $p$.
Also, since  ${\bf C}_G(a)\cap {\bf C}_G(b)\leq {\bf C}_G(ab)$, it follows
that $|G:{\bf C}_G(ab)|\leq p^2$.

Suppose that  
$|G:{\bf C}_G(ab)|= p^2$. Since $(ab)^G\subseteq a^G b^G$ and 
$|a^G|=|b^G|=p$, it follows then that $(ab)^G=a^G b^G$. Thus we may 
assume that $|G:{\bf C}_G(ab)|\leq p$ and so $|(ab)^G|\leq p$.

Let $Z$ be a normal subgroup of $G$ of order $p$ contained in ${\bf C}_G(a)$. 
Observe such subgroup exists since $|{\bf C}_G(a)|>p$, 
otherwise $G$ would have order $p^2$ and thus it would be abelian and therefore
$|a^G|=1$. 
Observe that 
$Z$ is a central subgroup of $G$ since $G$ is nilpotent, $Z$ is normal in $G$ and
$|Z|=p$.
Fix
$g\in G\setminus {\bf C}_G(a)$.

\begin{claim}\label{claim1}
If  $|(aZ)^{G/Z}|=1$ or $|(bZ)^{G/Z}|=1$  then either $\eta(a^G b^G)=1$ or
$\eta(a^G b^G)=p$.
\end{claim}
\begin{proof}
By symmetry, observe that it is enough to prove that if
$|(aZ)^{G/Z}|=1$ 
 then either $\eta(a^G b^G)=1$ or
$\eta(a^G b^G)=p$.

Observe that  $|(aZ)^{G/Z}|=1$ implies that $a^g=az$ 
for some $z\in Z\setminus \{1_G\}$, and therefore  
 $[a,G]=Z$. Similarly, if $|(bZ)^{G/Z}|=1$ then 
 $[b,G]=Z$. Observe that then $|(aZ)^{G/Z}|=|(bZ)^{G/Z}|=1$
 implies that $a^G b^G=\{ (ab)z\mid z\in Z\}$ and thus 
 either $[ab,G]=Z$, and so $a^G b^G=(ab)^G$, or
 $[ab,G]\cap Z=\{1_G\}$. If $[ab,G]\cap Z= \{1_G\}$, 
 by Lemma \ref{intersection} we get that $\eta(a^G b^G)=|Z|=p$.
 We may assume then that $|(bZ)^{G/Z}|=p$ and thus 
 $[b,G]\cap Z=\{1_G\}$. Since
 $a^g=az$,  given any integer $i$ we have that
 $a^{g^i}=az^i$.
 Thus
 \begin{equation}\label{equation1}
 (ab)^G=\{ab[b,g^i]z^i\mid i=0,\ldots, p-1\},
 \end{equation}
 and
 \begin{equation}\label{equation22}
 a^G b^G=\{ab[b,g^i]z^j\mid i,j=0,\ldots, p-1\}.
 \end{equation}
 Since $|(abZ)^{G/Z}|$ is either $1$ or $p$,
 either $[ab,G]=Z$,  and so $a^Gb^G=(ab)^G$,
 or $[ab,G]\cap Z=\{1_G\}$. If $[ab,G]\cap Z=\{1_G\}$, by \eqref{equation1},
  \eqref{equation22} and 
  Lemma \ref{intersection} we get that $\eta(a^G b^G)=p$. Thus
either 
$\eta(a^G b^G)=1$ or $\eta(a^G b^G)=p$. 
\end{proof}

By the previous claim, we may assume that $|(aZ)^{G/Z}|=|(bZ)^{G/Z}|=p$.
Since $|a^G|=|b^G|=p$, we get that
\begin{equation*}\label{abandz}
[a,G]\cap Z =\{1_G\}\mbox{ and } [b,G]\cap Z=\{1_G\}.
\end{equation*}
 By induction on the order of $G$, we have that 
  either $\eta((aZ)^{G/Z} (bZ)^{G/Z})=1$ or $\eta((aZ)^{G/Z} (bZ)^{G/Z})\geq \frac{p+1}{2}$.
  If $\eta((aZ)^{G/Z} (bZ)^{G/Z})\geq \frac{p+1}{2}$, by Lemma \ref{lemma1}
  we have that $\eta (a^G b^G)\geq \frac{p+1}{2}$. We may assume then that 
   $\eta((aZ)^{G/Z} (bZ)^{G/Z})=1$. 
   
   Observe that  $|[(ab)Z,G/Z]|> 1$ since 
    $\eta((aZ)^{G/Z} (bZ)^{G/Z})=1$ and $|(aZ)^{G/Z}|=|(bZ)^{G/Z}|=p$. Thus
    $|[(ab)Z,G/Z]|=p$. By Lemma \ref{samesize} we have that
    $[aZ,G/Z]=[bZ,G/Z]=[(ab)Z, G/Z]$ is a normal
   subgroup of $G/Z$.
   
    Let $Y$ be a normal subgroup of $G$ such that
   $Y/Z=[aZ, G/Z]$. Since  $|[a,G]|=|(aZ)^{G/Z}|=p$, 
   we have that $Y$ is a normal 
   subgroup of $G$ of order $p^2$ and $[a,G]Z=Y$. Set $[a,g]=y$.
   
   \begin{claim}\label{ycentral} If ${\bf C}_G(Y)\neq {\bf C}_G(a)$ or
   ${\bf C}_G(Y)\neq {\bf C}_G(b)$,
   then either $\eta(a^G b^G)=1$ or
$\eta(a^G b^G)= p$. 
  \end{claim}
  \begin{proof}
  Since $Y$ is normal in $G$, $Z$ is central and $|Y:Z|=p$, 
   either $|G:{\bf C}_G(Y)|=p$ or $Y$ is contained in the center ${\bf Z}(G)$ 
   of $G$.

Assume that $|G:{\bf C}_G(Y)|=p$ and ${\bf C}_G(Y)\neq {\bf C}_G(a)$.
 Fix $h\in {\bf C}_G(a)\setminus {\bf C}_G(Y)$. Observe such element
 exists since ${\bf C}_G(Y)\neq G$ and ${\bf C}_G(Y)\neq {\bf C}_G(a)$.
 Observe that $[h,y]\in Z\setminus\{1_G\}$. Set $[h,y]=z$.
 We can check that
 $[h^i,y]=z^i$ for $i=0,\ldots, p-1$. Since $[a,g]=y$, we have  
 $$(a^g)^{h^i}=(ay)^{h^i} =a^{h^i} y^{h^i}=ayz^i.$$ 
 We can check also that $(a^{g^i})^{h^j}$ for $i,j=1,2,\ldots, p$ are distinct elements.
 But then $|a^G|=p^2$.  Thus, if ${\bf C}_G(a)\neq{\bf C}_G(Y)$
 then ${\bf C}_G(Y)=G$.
 
 Observe that $[b,g]\in Y$ since $[bZ, G/Z]=Y/Z$. Let $[b,g]=y^rz^s$ for some integers
   $r, s\in \{0,\ldots, p-1\}$.
  We can check that then  $[a,G]=\{y^i\mid i=0,\ldots, p-1\}$,
$ [b,G]=\{ y^{ri} z^{si}\mid i=0,\ldots, p-1\}$,
$[ab,G]=\{ y^{(r+1)i} z^{si}\mid i=0,\ldots, p-1\}$ and
$[a,G][b,G]=\{ y^{ri+j} z^{si}\mid i,j=0,\ldots,p-1\}$. 
If $s\equiv 0\modu p$, then $[a,G][b,G]=[ab,G]$ and so $\eta(a^G b^G)=1$.
Otherwise by Lemma \ref{productsn} we have that
$a^G b^G= \cup_{v \in Z} (ab)^G v$,  and since $[ab,G]\cap Z=\{1_G\}$,
by Lemma \ref{intersection} it follows that $\eta(a^G b^G)=p$.

Similarly, we can check that 
if ${\bf C}_G(Y)\neq {\bf C}_G(b)$, then either $\eta(a^G b^G)=1$ or
$\eta(a^G b^G)=p$.
  \end{proof}
  
  By the previous claim, we may assume that ${\bf C}_G(a)={\bf C}_G(Y)={\bf C}_G(b)$.
  
   Set $[y,g]=z$. Since $g\in G\setminus {\bf C}_G(a)$, and
    by the previous claim 
   it follows that $z\in Z\setminus \{ 1_G\}$. Because $Z$ is central, $[y,g^i]=z^i$ for 
   any integer $i$. Since $[a,g]=y$, we can check by
   induction on $i$ that
   \begin{equation}\label{agi}
   a^{g^i}= ay^i z^{\frac{i(i-1)}{2}} \mbox{ for any integer }i\geq 0.
   \end{equation} 
   
  Since $a^{g^p}=a$ and $z^p=1_G$, it follows that $y^p=1_G$. 
  Thus $Y$ is an elementary abelian
  normal 
  subgroup of $G$. 
    
    Let $b^g=by^rz^t$ for some integers $r$ and $t$.
     Since $|b^G|=p$ and $[b,G]\cap Z=\{1_G\}$, we have that
    $r\not\equiv 0 \modu p$.  We can check by induction on $i$ that 
   \begin{equation}\label{bgi}
   b^{g^i}=by^{ri} z^{r\frac{i(i-1)}{2} + ti} \mbox{ for any integer }i\geq 0.
    \end{equation}
    
    Since $b\in {\bf C}_G(b)={\bf C}_G(a)={\bf C}_G(Y)$,
     combining \eqref{agi} and \eqref{bgi}, we get 
   \begin{equation}\label{212}
   (ab)^{g^i}= (ab)y^{(r+1)i} z^{(r+1)\frac{i(i-1)}{2}+ti}.
   \end{equation} 
   Since $[(ab)Z, G/Z]=[aZ,G/Z]$  and $|[aZ,G/Z]|=p$, we have that
   $[ab,G]\cap Z=\{1_G\}$. Thus we must have that $r+1\not\equiv 0\modu p$.
    Fix $i\in \{0,\ldots, p-1\}$.    Let $j\in \{0,\ldots, p-1\}$
    such that $i\equiv (r+1)j\modu p$. Then by \eqref{212} we have 
    \begin{equation}\label{counting1}
    ab^{g^i}= aby^iz^{\frac{i(i-1)}{2}}= (ab)^{g^j} z^{s_i},
    \end{equation}
    \noindent where $s_i= \frac{i(i-1)}{2}-(r+1)\frac{j(j-1)}{2}-tj$. Set
    $S=\{s_i\mid i=0, \ldots, p-1\}$. 
    Observe that 
    \begin{eqnarray*}\label{counting2} 
    s_i&=&\frac{i(i-1)}{2}-(r+1)\frac{j(j-1)}{2}-tj\\
    &=&\frac{(r+1)j((r+1)j-1)}{2}-(r+1)\frac{j(j-1)}{2}-tj\\
    &=&\frac{(r+1)r}{2}j^2-tj. 
    \end{eqnarray*}
    Since $r+1\not\equiv 0\modu p$ and $i\equiv (r+1)j\modu p$, 
     $S\modu p=\{s_i\modu p\mid i=0,\ldots ,p-1\}=
    \{\frac{r(r+1)}{2}j^2-tj\modu p\mid j=0,\ldots, p-1\}$.
     Therefore by Lemma \ref{counting}
    we have that either $S\modu p$ has one element, or $S\modu p$
     has at least $\frac{p+1}{2}$. By \eqref{counting1}, we have
 \begin{equation}\label{counting3}
  a^G b^G= \cup_{i=0}^{p-1}(ab^{g^i})^G= \cup_{i=0}^{p-1}(ab)^G z^{s_i}.
  \end{equation}
   Since $Z=<z>$ is a central subgroup of $G$ of order $p$, $[ab,G]\cap Z=\{1_G\}$,
  and the set $S\modu p$ has either one element or at least
  $\frac{p+1}{2}$ elements, by Lemma \ref{intersection}
   and \eqref{counting3} we have that either 
   $\eta(a^G b^G)=1$ or  $\eta(a^G b^G)\geq \frac{p+1}{2}$ and the proof of
   the theorem is
   complete.    
\end{proof}
\end{section}

\begin{section}{Proof of Theorem B}
The following two results are  proved in \cite{edithhomogen}.
\begin{lem}\label{eta1prop}
Let $G$ be a finite group, $a^G$ and $b^G$ be  conjugacy classes of $G$. Assume that 
${\bf C}_G(a)={\bf C}_G(b)$.
 Then  $a^G b^G=(ab)^G$ if and only if  $[ab,G]=[a,G]=[b,G]$ and $[ab,G]$ is a normal 
subgroup of $G$. 
   \end{lem}
\begin{lem}\label{notinthecenter}
Let $G$ be a group of odd order and $a^G$ be the conjugacy class of $a$ in $G$.
Then 
\begin{equation}\label{centerg}
{\bf Z}(G)\cap a^G a^G\neq \emptyset
\end{equation}
\noindent if and only if $|a^G|=1$.
Thus if $|a^G|>1$ and  $b^G\subseteq a^Ga^G$, then $|b^G|>1$.
\end{lem}
 \begin{proof}[Proof of Theorem B]
 Since $p$ is a prime number and $G$ is a nilpotent group, by Lemma \ref{observation1}
 we may assume that $G$ is a $p$-group.
 
 Since $|a^G|=p$ and $|a^G|=|G:{\bf C}_G(a)|$,  ${\bf C}_G(a)$ is a normal subgroup of 
 index $p$. Clearly $a\in {\bf C}_G(a)$, and since ${\bf C}_G(a)$ is normal in $G$, 
 $a^g\in {\bf C}_G(a)$ for any $g\in G$. Thus 
 $aa^g=a^ga$ for any $g\in G$.
 
  Fix $g\in G\setminus {\bf C}_G(a)$.
  Since $aa^{g^i}$ is $G$-conjugate to
  $a^{g^{p-i}} a$, $aa^{g^{p-i}}= a^{g^{p-i}}a$ and 
  $a^G a^G= \cup_{i=0}^{p-1} (aa^{g^i})^G$, 
  it follows that $\eta(a^G a^G)\leq \frac{p+1}{2}$.
  Thus by Theorem A we have that either 
  $\eta(a^G a^G)=1$ or $\eta(a^G a^G)=\frac{p+1}{2}$. 
  
  If $\eta(a^G a^G)=1$, then by Lemma \ref{eta1prop} we have that
  $[a,G]=[a^2,G]$ is a normal subgroup of $G$ and i) holds. 
  
  We may assume now that $\eta(a^G a^G)=\frac{p+1}{2}$. Since $G$ is of odd order, 
  by Lemma \ref{notinthecenter}, we have that $a^Ga^G$ must be the union 
  of conjugacy classes of size at least $p$.
  Since $|a^G a^G|\leq p^2$ and $\eta(a^G a^G)=\frac{p+1}{2}$,
  it follows that $a^G a^G$ is the union of $\frac{p+1}{2}$ distinct
  conjugacy classes of size $p$ and ii) holds.
 \end{proof}
 \end{section}
\begin{section}{Examples}

\begin{hypo}\label{basic}
Let $p$ be an odd prime and $n>0$ be an integer.
Let $G_0$ be a $p$-group and $g_0^{G_0}$ be a conjugacy class of $G$ with 
$g_0^{G_0} g_0^{G_0}=(g_0^2)^{G_0}$ and $|g_0^{G_0}|=p^{n-1}$. 

Let $N=G_0\times G_0\times \cdots\times G_0$ be the direct product of $p$ copies of $G_0$.
Set $$a=(g_0, e,\ldots, e),$$ 
\noindent where $e$ is the identity of $G_0$.
Let
$C=<c>$ be a cyclic group of order $p$. 
 Observe that $C$ acts on $N$ by
 \begin{equation}\label{actionc}
 c:(n_0, n_1,\ldots,n_{p-1})\mapsto(n_{p-1},n_0, \dots, n_{p-2})
 \end{equation}
 \noindent for any $(n_0, n_1, \ldots,n_{p-1})\in N$.
 
Let $G$ be the semi direct product of $N$ and $C$, i.e $G$ is the wreath product of $G_0$ and 
$C$. 
\end{hypo}
{\bf Remark.} Fix a prime $p>2$ and an integer $l>0$.
Let $E$ be an extraspecial group of exponent $p$ and order $p^{2l+1}$. Let $g\in E\setminus {\bf Z}(E)$. We can check that $g^E=g{\bf Z}(E)$ and thus 
$g^E g^E=(g^2)^E$. 
Clearly if $G$ is an abelian group and $g\in G$, then $g^Gg^G=(g^2)^G$.
Therefore given any prime $p$ and any integer $l\geq  0$, there exist a $p$-group
$G$ and a conjugacy class $a^G$ such that $a^G a^G=(a^2)^G$ and
$|a^G|=p^l$, and so there always exists a group satisfying the previous hypothesis.

 \begin{prop}\label{basico}
Assume Hypotheses \ref{basic}. 
 Then $\eta(a^G a^G)=\frac{p+1}{2}$ and $|a^G|=p^n$. Thus given any 
integer $n>0$ and any prime $p>2$, there exist some $p$-subgroup $G$ and
some conjugacy class $a^G$ of $G$ such that $\eta(a^G a^G)=\frac{p+1}{2}$ and
$|a^G|=p^n$. 
\end{prop}
\begin{proof}
Observe that $a^G=\{(x, e, \ldots,e), (e,x,\ldots e),\ldots, (e,\ldots, e,x)\mid x\in \ g_0^{G_0}\}$. Thus $|a^G|=p|g_0^{G_0}|=p(p^{n-1})=p^n$.

Observe also that 
\begin{equation*}
\begin{split}
a^Ga^G=& \{(xy, e,e, \ldots, e)^b,(x,y,e,\ldots,e)^b, (x,e,y,e,\ldots, e)^b,\ldots\\
&(x,e,\ldots, e, y)^b\mid x, y \in g_0^{G_0}\mbox{ and } b\in C\}.
\end{split}
\end{equation*}
Since $g_0^{G_0} g_0^{G_0}=(g_0^2)^{G_0}$, we have that
\begin{equation}\label{g0square}
  (a^2)^G=\{(x^2, e,\ldots,e)^b \mid x \in g_0^{G_0}\mbox{ and } b\in C\}.
\end{equation} 

Set $a_i=(g_0,e,\ldots, g_0, \ldots,e)$, where the second $g_0$ is in the 
$i$-position, for example $a_0=(g_0 g_0, e,\ldots,e)$ and $a_1=(g_0,g_0,e,\ldots, e)$. 
Now observe that $(x,y, e,\ldots,e)$ is $G$-conjugate to 
$(u,e, \ldots, e, v)$ for any $x, y , u, v\in g_0^{G_0}$. Furthermore, we can check that
$(x,e, \ldots, e,y,e \ldots, e)$,
where $y$ is in the $i$-th position, is $G$-conjugate to
$(u, e,\ldots, v, e)$, where $v$ is in the $j$-position if and only if
$j\equiv\pm i \modu p$. Thus $a_i^G\cap a_j^G\neq \emptyset$ if and only
if $j\equiv \pm i\,\modu p$. 
We can check that 
\begin{equation*}
a^G a^G=\cup_{i=0}^{\frac{p-1}{2}} a_i^G.
\end{equation*}
Therefore $\eta(a^G a^G)=1+ \frac{p-1}{2}=\frac{p+1}{2}$.
\end{proof}
{\bf Remark.} 
Assume Hypotheses \ref{basic}. Let 
$b=(g_0,g_0, e,\ldots ,e)\in N$.
We can check that $\eta(a^G b^G)=p-1$.
The author wonders if there exist a prime $p>5$, a $p$-group $G$, and
conjugacy 
classes $a^G$ and $b^G$ of $G$ with $|a^G|=|b^G|=p$ such that
 $\frac{p+1}{2}< \eta(a^G b^G)< p-1$.

\begin{hypo}\label{hypot}	Fix a prime $p$ and let
$F=\{0,1,\ldots,p-1\}$ be the finite field with $p$ elements. Observe that $F$ is also
a vector
space of dimension 1 over itself.
Let 
$ A=\Aff(F)$ be the affine group of $F$.  
Observe that the group $A$ is a cyclic by cyclic group and thus
it is supersolvable. 

Let $C$ be a cyclic group of order $p$. 
 Set $X=F$ and  $K = C^X = \{ f: X \to C \}$. Observe that
$K$ 
is a group via pointwise multiplication, and clearly
$A$ acts on this group (via its action on $X$).
 
 	Let $G$ be the wreath product of $C$ and $A$ relative to $X$,
i.e. $G = K \rtimes A$. We can check  that $G$ is a supersolvable group.
\end{hypo}
\begin{prop}\label{super}
Assume Hypotheses \ref{hypot}. Let $e$ be the identity of $C$.
Set $a=(c,e,e,\ldots, ,e)$ in $K$, where $c\neq e$. 
Then $a\in G$, the conjugacy class $a^G$ of $\,G$  has size $p$ and 
$a^G a^G = (a^2)^G \cup (c,c,e,  \ldots, e)^G$. 
Thus  $\eta(a^G a^G)=2$.

Therefore, given any prime $p$, there exist a supersolvable group $G$  and
a conjugacy class $a^G$ of $G$ with $|a^G|=p$ and $\eta(a^G)=2$.
\end{prop}
\begin{proof}
Observe that 
$a^G=\{(c, e, \ldots,e), (e,c,e,\ldots, e),\ldots, (e,e,\ldots,e,c)\}$. 
Thus $a^G$ has $p$ elements.
 We can check that
\begin{equation}\label{ag21}
(a^2)^G=\{(c^2, e, \ldots,e), (e,c^2,e,\ldots, e),\ldots, (e,e,\ldots,e,c^2)\}.
\end{equation}
\noindent Also, we can check that
the set 
$(c, c, e, \ldots , e)^G$ is the set of all elements of the form
 $(e,\ldots, e, c,e,  \ldots, e,c,e,\ldots,e)$,
 where the 
two $c$'s  are in the $i$-th and the $j$-th position for some $i$, $j$ with $1\leq i<j\leq p$.
Observe that 
\begin{equation*}\label{ag2}
a^G a^G  = \{(c^2, e, \ldots,e), (e,c^2,e,\ldots, e),\ldots, (e,e,\ldots,e,c^2)\} \cup (c, c, e, \ldots , e)^G .
\end{equation*}
\noindent 
 Thus  
$a^G a^G= (a^2)^G \cup (c, c, e, \ldots , e)^G$.
Clearly $(a^2)^G \neq (c, c, \ldots , e)^G$. We conclude then that
$\eta(a^G a^G)=2$.  
\end{proof} 
 \end{section}

\end{document}